# Study of Air Pollution Impact on Human Health in Major Cities in India using Fractal Analysis

Santanu Nandi

Dept. of Mathematics, School of Advanced Sciences, VIT - AP University, Andhra Pradesh - 522241, India.

Email: santanu282@gmail.com

**ABSTRACT:**

This study aims to examine the historical air pollution data from major Indian cities using fractal analysis to measure environmental risk. The fractal dimension of the major air pollutants is computed to evaluate the volatility and complexity of air quality patterns in Delhi, Kolkata, Mumbai, and Bengaluru. Fractal statistics is applied to analyze the fractal spectrum for each region, and further divided into different scales to capture the variation in air quality behaviour. This analysis reveals distinct fractal patterns for each city, showing different levels of environmental instability. These findings make fractal analysis a valuable tool for comparing air quality dynamics. This study offers a new way to measure changes in air quality, helping policymakers create targeted solutions for each city.

## 1. INTRODUCTION:

Air pollution is an increasingly critical environmental issue, especially in rapidly urbanizing regions where high levels of particulate matter (PM) and other pollutants pose a significant threat to human health and ecological balance. Urban centres such as Delhi, Mumbai, Kolkata, and Bengaluru face severe pollution challenges, making air quality assessment and management essential to public health and environmental sustainability. Traditional methods for analyzing air quality often rely on standard statistical tools that assume normal distribution. However, these methods may overlook the complex, non-linear, and scale-invariant patterns found in real-world environmental data [1][4][7].

Fractal analysis offers an innovative approach to understanding air quality fluctuations. Fractal geometry, introduced by Mandelbrot, has proven to be effective in capturing irregular, self-similar patterns within data. The use of fractals provides insights into data characterized by scale invariance and fractional dimensions, both of which are prevalent in environmental datasets [1][2]. In recent years, fractal techniques have been applied to assess volatility in various domains, from stock market dynamics to air pollution, illustrating their versatility in handling data that deviate from normal distribution assumptions [2][3][4].

This research leverages fractal analysis to investigate AQI variability across different cities, employing techniques such as fractal dimension and multifractal scaling to capture the inherent complexity and risk levels associated with air quality. By analyzing fractal patterns in AQI data, this study aims to reveal the temporal and spatial dynamics of pollution, identify scale-dependent patterns, and highlight potential health risks across these cities.

# 2. PRELIMINARIES:

Fractal geometry, a branch of mathematics dealing with irregular and self-similar patterns, has evolved to include fractal statistics, which allows for the study of complex datasets exhibiting non-normal characteristics. According to Padua et al., fractal geometry's features-self-similarity, scale invariance, and fractal dimension-provide a foundation for understanding irregular datasets, which are typical in fields such as finance, meteorology, and environmental science [1]. The development of fractal statistics enables researchers to analyze datasets with non-existent means and variances, properties that traditional methods cannot handle effectively [1][2].

The application of fractals in finance has demonstrated the relevance of fractal statistics for understanding market risks. Studies examining the volatility of stock prices have shown that fractal dimensions can reflect the risk associated with different stocks, with higher fractal dimensions indicating greater volatility [2]. This approach has proven useful in identifying patterns in seemingly chaotic data, allowing analysts to evaluate risk based on the degree of fluctuation or "ruggedness" in stock prices. Similar techniques have been extended to environmental datasets, where fractal dimensions offer a measure of the variability in pollution data, reflecting the potential health risks associated with high volatility in air quality [2][3].

Fractal methods are increasingly used to study air pollution patterns, with research showing that pollutants like PM2.5 and PM10 often exhibit complex, multifractal behaviours in urban settings [4][5][6]. Studies in cities like Shanghai and in regions like Kazakhstan have highlighted how fractal analysis captures pollutant variability that traditional methods might miss, offering insights into persistent and scale-dependent patterns in air quality data [5]. In India, fractal analysis has proven valuable for AQI studies, particularly for revealing seasonal and daily variations in pollutant behaviour. For example, fractal properties of PM2.5 across the Indo-Gangetic Plain shift with weather patterns, while studies in New Delhi show fractal dimensions can indicate periods of high pollution risk by capturing the interactions among major pollutants [3][4][8].

Research further demonstrates that pollutant data frequently show multifractal patterns, which indicate changes in pollutant concentration across scales. These findings suggest that fractal and multifractal analysis can improve pollution forecasting and risk assessment by accounting for variability in pollutant behaviour across time and location [6][7].

# 3. METHODOLOGY:

### *3.1. Dataset Overview:*

The AQI data utilized in this study was sourced from the Central Pollution Control Board (CPCB), the Indian government body responsible for air quality monitoring across the nation. CPCB manages a network of air quality monitoring stations operated in collaboration with central and state-level Pollution Control Boards as well as local agencies, providing a comprehensive and coordinated approach to air quality assessment.

The dataset spans from 2017 to September 2024, offering a longitudinal view of daily AQI trends. Data was collected from four major Indian cities—Delhi, Mumbai, Kolkata, and Bengaluru—each representing diverse environmental and geographic conditions. Each dataset includes two columns: a time-stamp indicating daily observation dates, and the AQI value, reflecting the average pollution level for each city on a daily basis (Fig.1). The data was accessed through the CPCB's official portal, ensuring consistency and reliability in the information across all years and locations. (Dataset Page Link)

| Delhi | Mumbai | Kolkata | Bengaluru |
|---|---|---|---|
| Time-stamp  AQI | Time-stamp  AQI | Time-stamp  AQI | Time-stamp  AQI |
| 0  2017-01-01  345.0 | 0  2017-01-01  276.0 | 0  2018-04-01  60.0 | 0  2017-01-01  84.0 |
| 1  2017-01-02  337.0 | 1  2017-01-02  171.0 | 1  2018-04-02  60.0 | 1  2017-01-02  74.0 |
| 2  2017-01-03  331.0 | 2  2017-01-03  192.0 | 2  2018-04-03  60.0 | 2  2017-01-03  69.0 |
| 3  2017-01-04  381.0 | 3  2017-01-04  177.0 | 3  2018-04-04  60.0 | 3  2017-01-04  71.0 |
| 4  2017-01-05  321.0 | 4  2017-01-05  186.0 | 4  2018-04-05  60.0 | 4  2017-01-05  74.0 |

**Figure 1: First 5 observations of each city**

## 3.2. Data Preprocessing:

To prepare the AQI dataset for analysis, various preprocessing steps were implemented to address missing values and ensure data consistency. The dataset contained some missing AQI values, which were handled using median imputation. The median was selected for its robustness against extreme values, making it particularly suitable for historical time series data such as AQI. When handling missing values, an approach based on the quantity and distribution of missing data was followed. For example, if a month had relatively few missing values, these gaps were filled with the median AQI for that specific month to preserve seasonal characteristics. For larger gaps, imputation was performed using the median AQI for the corresponding year. In cases where a significant amount of data was missing, these observations were removed to avoid distorting the dataset's natural distribution.

Outlier detection was conducted to examine the dataset for unusually high or low AQI values. Minimum and maximum AQI values were reviewed to ensure they reflected realistic fluctuations in air quality. No significant outliers were found, and occasional spikes and dips in AQI were attributed to natural factors such as weather changes, traffic density, and industrial activity. As such, these fluctuations were retained in the dataset to accurately reflect the variability inherent in air quality data.

## 3.3. Analysis of Data:
*Indicators of Fractal Distribution*

To justify the use of fractal statistics for analyzing AQI data, it was essential to confirm that the datasets did not follow a normal distribution. Fractal analysis is best suited to datasets characterized by irregular, non-linear, and non-normal patterns, as it can reveal complex dynamics that standard statistical methods may not capture. By establishing the non-normality of the AQI datasets, the groundwork for fractal analysis was effectively laid.

As an initial step, histograms were plotted for each city's AQI dataset to visually assess distribution patterns in the data (Fig.2.).

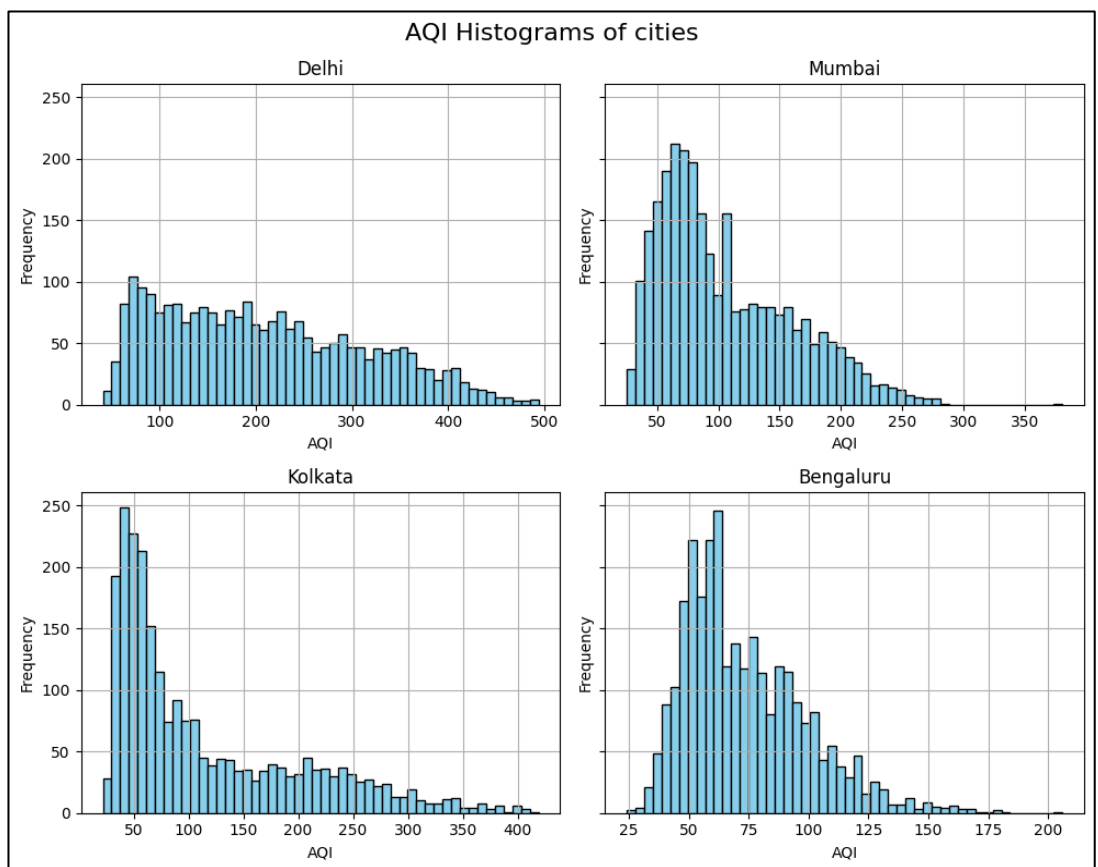

**Figure 2: AQI Histograms of cities**

These histograms provide an overview of the AQI values' distribution across cities. Unlike the bell-shaped curve expected in a normal distribution, the AQI histograms exhibit clear right-skewness in each, suggesting the absence of normal distribution characteristics. This first analysis through histograms revealed that the AQI data's distribution does not align with the assumptions of normality, which hints at underlying complexity that fractal analysis can further explore.

Following the histogram inspection, Q-Q (Quantile-Quantile) plots were generated for each city's dataset to provide a more detailed visual comparison between observed AQI values and those expected under a normal distribution (Fig.3.).

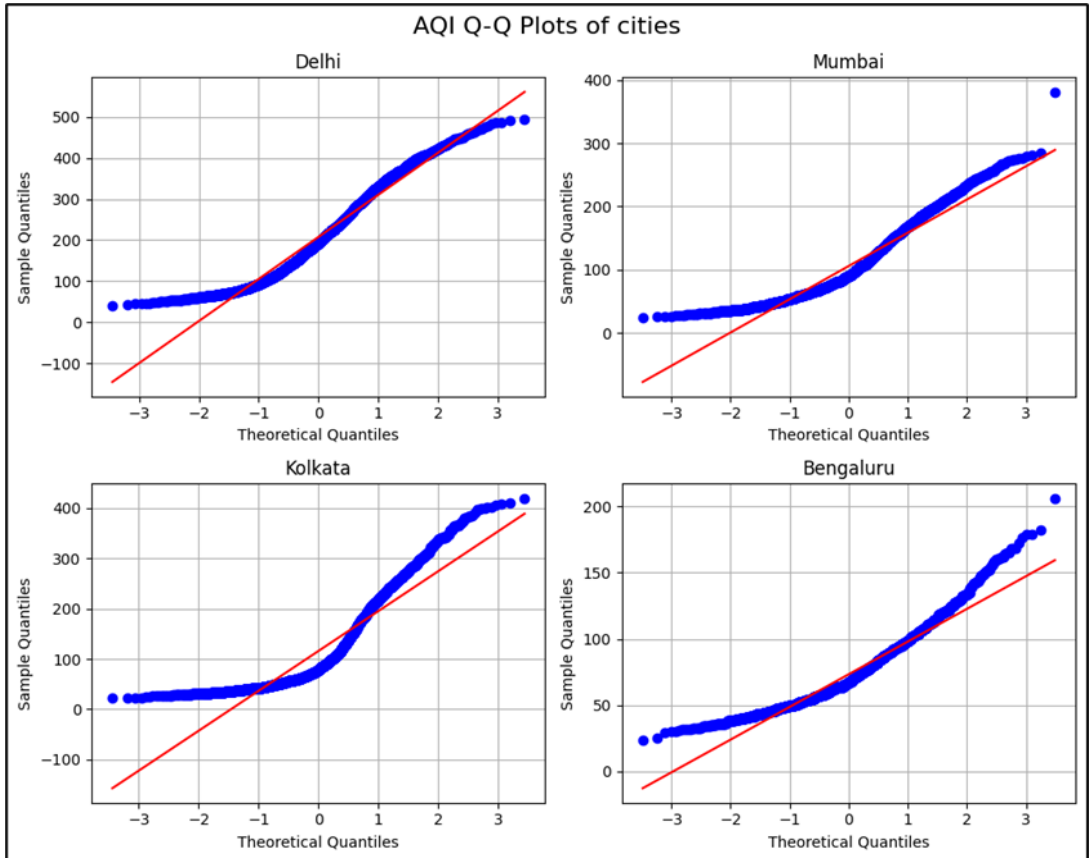

**Figure 3: AQI Q-Q Plots of cities**

The Q-Q plots further underscore the findings from the histograms. Deviations from the expected normal line are evident, as the AQI data points in each plot stray significantly from the diagonal, showing that the data distributions diverge from a normal distribution model. This deviation from the diagonal line, combined with the skewed histograms, confirms that the AQI values do not fit the normal distribution model.

To validate these observations statistically, the Kolmogorov-Smirnov (K-S) Goodness-of-Fit Test was applied to each AQI dataset (Table 1).

**Table 1: Kolmogorov-Smirnov Goodness of Fit Test Results**

| AQI City | Mean | Std. Dev. | N | K-S Statistic | P-Value |
|---|---|---|---|---|---|
| Delhi | 207.764 | 104.394 | 2465 | 0.073 | 4.65e-12 |
| Mumbai | 106.109 | 54.522 | 2830 | 0.118 | 5.37e-35 |
| Kolkata | 115.764 | 86.038 | 2375 | 0.173 | 4.02e-63 |
| Bengaluru | 73.386 | 25.416 | 2830 | 0.107 | 4.65e-29 |

The K-S test results provided additional confirmation, with p-values very close to 0, indicating strong evidence against the normality of each AQI dataset. This combination of histogram and Q-Q plot analyses, along with the statistical support from the K-S test, strongly validates the non-normality of the datasets. Consequently, fractal statistics can be applied to the AQI datasets to further analyze their scale-dependent, irregular patterns.

## 3.4. Fractal Model:
*Applying Fractal Statistics*

To capture the complex, irregular patterns of AQI data, fractal statistics offer an ideal approach by quantifying the inherent ruggedness and variability within the data. Unlike conventional statistical methods, which assume normal distribution and regularity, fractal models allow for the analysis of datasets that display self-similarity, scale invariance, and unpredictable fluctuations. This study applies a fractal framework to analyze the AQI time series data from major Indian cities, as traditional statistical tools may overlook the intricate patterns present in environmental data.

The process involves several key steps: calculating the fractal dimension, which serves as a quantitative measure of the AQI data's complexity; implementing an algorithm to systematically compute fractal properties across the datasets; and generating the fractal spectrum to visualize scale-dependent behaviours. Additionally, a scale disparity analysis is conducted to compare variations across small, medium, and large scales, allowing for deeper insights into the relative environmental risks associated with air quality in different regions. This application of fractal statistics thus provides a comprehensive view of the underlying patterns in AQI data, enabling a more accurate understanding of air quality dynamics.

### 3.4.1. Calculating Fractal Dimension

To analyze the AQI data's complexity and variability across different cities, fractal dimensions were computed as a key indicator of environmental riskiness. Using the approach proposed by Padua et al. (2013), the fractal dimension ($\lambda$) was calculated for each AQI dataset, providing insight into the inherent ruggedness and scale-invariance of air quality patterns. The values in each dataset were first tabulated and ranked in ascending order, creating a structured basis for fractal analysis.

The fractal dimension ($\lambda$) was calculated using the following formula:

$$\lambda = 1 - \frac{log(1-\alpha)}{log\left(\frac{x_i}{\theta}\right)}$$

$where: \alpha = \frac{rank(x_i)}{n}$

$and: \theta = \min(x_1, x_2, \ldots, x_n)$

Following this methodology, the minimum AQI value ($\theta$) was identified in each dataset to serve as a reference point, allowing for a comparative view of scaling behaviours. This value acted as a "microscope" that sharpens the view of finer details in the AQI data, highlighting both the variability and stability of air quality across the analyzed cities. By using a uniform $\theta$ value, a consistent comparative framework was achieved, ensuring that any observed differences in fractal dimensions reflect true variations in AQI patterns rather than data-specific artifacts.

*Algorithm:*

To quantify the complexity of AQI patterns across the datasets, the following algorithm was employed to calculate the average fractal dimension $\bar{\lambda}$, which provides a measure of environmental variability for each city.

Input: AQI dataset $X = (x_1, x_2, \ldots, x_n)$

1. Identify Minimum Value ($\theta$):
    - Determine the minimum AQI value ($\theta$) in the dataset to serve as a baseline for scaling.

2. Sort AQI Values:

    o   Arrange the AQI values in ascending order.

3. Calculate Fractal Dimension ($\lambda$) for Each Observation:

    o   For each AQI value, assign it a rank based on its position in the sorted list.

    o   Compute the fractal dimension ($\lambda$) corresponding to each ranked AQI value.

4. Average Fractal Dimension ($\bar{\lambda}$):

    o   Calculate the average of all $\lambda$ values to obtain the overall fractal dimension for the dataset.

Output: The average fractal dimension $\bar{\lambda}$ provides a summary of the dataset's fractal characteristics.

The following table provides a summary of descriptive statistics for the fractal dimension ($\lambda$) values across four cities, offering insight into the complexity and irregularity in their AQI patterns.

**Table 2: Descriptive Statistics of Fractal Dimension ($\lambda$)**

| AQI City | N | Mean ($\bar{\lambda}$) | Std. Dev. | Minimum | Median | Maximum |
|---|---|---|---|---|---|---|
| Delhi | 2463 | 1.546 | 0.394 | 1.013 | 1.448 | 4.145 |
| Mumbai | 2828 | 1.609 | 0.412 | 1.016 | 1.540 | 4.266 |
| Kolkata | 2373 | 1.562 | 0.318 | 1.018 | 1.546 | 3.657 |
| Bengaluru | 2828 | 1.773 | 0.555 | 1.006 | 1.673 | 4.923 |

As shown in Table 2, the mean fractal dimension ($\bar{\lambda}$) varies by city, with Delhi showing $\bar{\lambda}$ of 1.546 and a moderate range from 1.013 to 4.145, suggesting considerable irregularity in its AQI data. Mumbai, with a slightly higher $\bar{\lambda}$ of 1.609 and a range up to 4.266, indicates a more complex AQI pattern than Delhi, suggesting a tendency toward higher variability. Kolkata, with an $\bar{\lambda}$ of 1.562 and a narrower range, shows somewhat less fluctuation in its AQI patterns, implying a relatively stable yet complex behavior in air quality. Bengaluru, however, stands out with the highest $\bar{\lambda}$ of 1.773 and the widest range (1.006 to 4.923), reflecting a pronounced level of irregularity and frequent fluctuations in AQI. This high fractal dimension and wide range in Bengaluru suggest that the city experiences more complex and varied air quality patterns, possibly influenced by a combination of diverse pollution sources and environmental conditions. These differences highlight the variability in AQI complexity across cities, warranting further exploration through the fractal spectrum to examine scale-dependent behaviours more closely.

*3.4.2. Fractal Spectrum*

Fractal observations are characterized by scale invariance and self-similarity, properties that allow for the analysis of data across different levels of detail. After examining the ruggedness of the AQI data through the fractal dimension ($\lambda$), establishing scale invariance is essential to fully understanding the patterns in air quality across cities. Scale (S) is defined as:

$$S = \frac{1}{log\left(\frac{x_i}{\theta}\right)}$$

where $x_i$ is the AQI observation and $\theta$ is the minimum AQI value in the dataset. This measure captures how variations in AQI behave over different scales, providing insight into the stability or variability of the data.

With the fractal dimensions and scale measures established, the next step involved constructing the fractal spectrum by plotting the relationship between scale and fractal dimensions ($\lambda$) for each city's AQI dataset. This spectrum illustrates how AQI values fluctuate across different scales, highlighting distinct behaviors in the air quality patterns of each city.

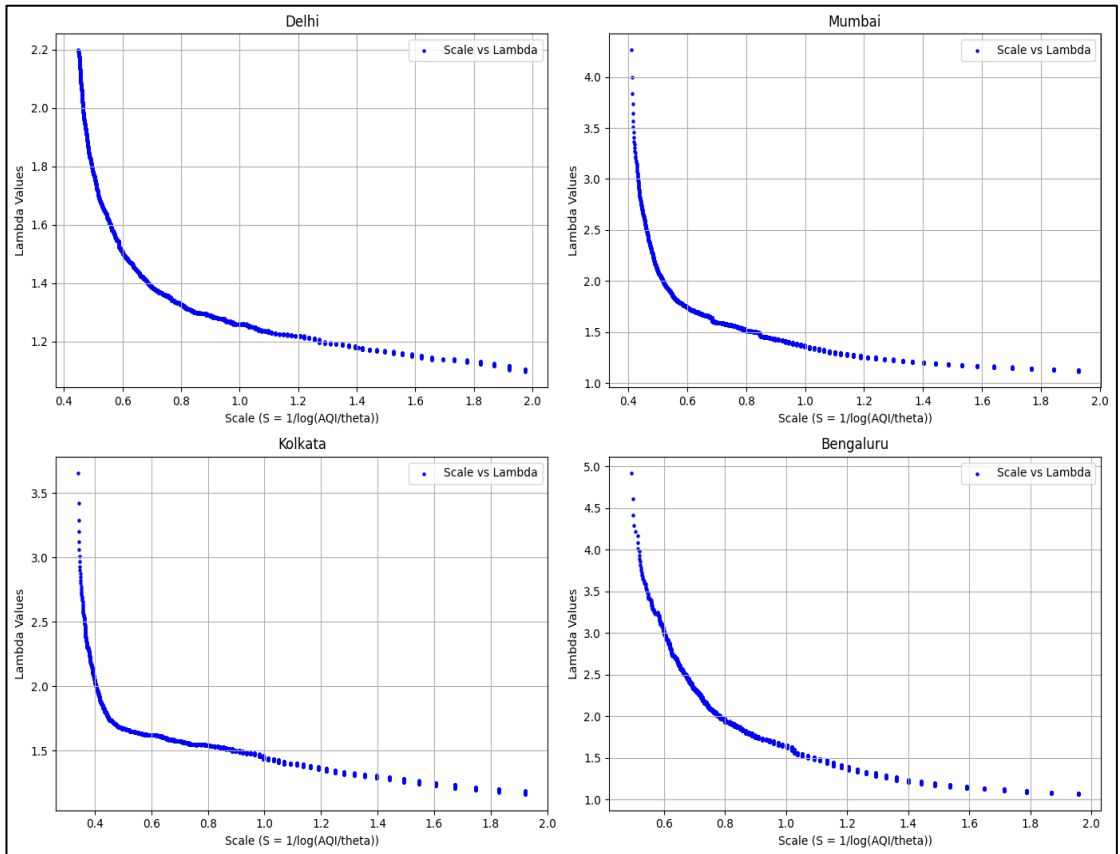

**Figure 4: Fractal Spectrum Plot of cities**

The fractal spectrum plots in Figure 4 highlight distinctive fractal behaviour in the AQI data across Delhi, Mumbai, Kolkata, and Bengaluru, consistent with the earlier analysis of fractal dimensions. Delhi and Kolkata exhibit relatively smoother decay curves, indicating a more stable air quality pattern with moderate volatility. This smoother pattern aligns with their lower average fractal dimensions, suggesting that while there are fluctuations, they are not excessively irregular.

In contrast, Bengaluru and Mumbai display more pronounced irregularities in their fractal spectra. Bengaluru, in particular, shows the highest degree of variation in the decay curve, reflecting its highest mean fractal dimension and a broad range of λ values, indicative of frequent and substantial AQI fluctuations. This pattern implies a high level of AQI instability and increased environmental risk. Mumbai also shows a moderately steep and uneven decay in its fractal spectrum, aligning with its slightly higher fractal dimension compared to Delhi and Kolkata, suggesting greater volatility but to a lesser degree than Bengaluru.

Overall, these fractal spectrum plots support the earlier findings: cities with smoother decay patterns (like Delhi and Kolkata) have relatively stable AQI, while cities with more irregular spectra (like Bengaluru and Mumbai) experience greater variability, suggesting higher air quality risks in these areas.

### 3.4.3. Scale Disparity Analysis

To further analyze the variability within AQI data, the scale disparity analysis provides insights into the ruggedness or volatility of air quality across different scales. For this purpose, the AQI data for each city was divided into three scales -small, medium, and large- based on AQI values. This division allowed for a closer examination of how AQI volatility changes across different levels. By analyzing the fractal dimension ($\lambda$) within these scales, we assess the degree of AQI fluctuation at each level.

In this analysis, the average AQI value in both the small and large scales was raised to the power of the average $\lambda$ in each respective scale. The difference between these values indicates the scale disparity, which provides an understanding of the degree of fluctuation present at each level of air quality. A larger disparity between small and large scales is indicative of higher volatility, suggesting an environment where air quality shifts significantly over time, potentially posing a greater risk to public health.

Table 3A: Scale Disparity Analysis, Delhi

| Delhi AQI | Average of data ($a$) | Average of lambda ($b$) | Ruggedness ($a^b$) |
|---|---|---|---|
| Small Scale | 317.555 | 1.854 | 43488.291 |
| Large Scale | 83.634 | 1.167 | 175.154 |
| Disparity | | | 43313.137 |

Table 3B: Scale Disparity Analysis, Mumbai

| Mumbai AQI | Average of data ($a$) | Average of lambda ($b$) | Ruggedness ($a^b$) |
|---|---|---|---|
| Small Scale | 185.315 | 2.169 | 83003.811 |
| Large Scale | 49.239 | 1.186 | 101.643 |
| Disparity | | | 82902.168 |

Table 3C: Scale Disparity Analysis, Kolkata

| Kolkata AQI | Average of data ($a$) | Average of lambda ($b$) | Ruggedness ($a^b$) |
|---|---|---|---|
| Small Scale | 245.309 | 1.966 | 49908.751 |
| Large Scale | 40.019 | 1.217 | 89.116 |
| Disparity | | | 49819.635 |

Table 3D: Scale Disparity Analysis, Bengaluru

| Bengaluru AQI | Average of data ($a$) | Average of lambda ($b$) | Ruggedness ($a^b$) |
|---|---|---|---|
| Small Scale | 108.476 | 2.533 | 143054.053 |
| Large Scale | 46.666 | 1.196 | 99.112 |
| Disparity | | | 142954.941 |

As seen in Table 3A to 3C, the scale disparity analysis provides a comparative measure of risk across the cities by evaluating the ruggedness difference between small and large scales. Delhi, Mumbai, Kolkata, and Bengaluru each show varying levels of disparity, which reflect differences in AQI variability and stability. Consistent with the findings from the fractal dimension calculations and fractal spectrum plots, Bengaluru demonstrates the highest disparity between scales, with a ruggedness difference of approximately 142,955. This substantial disparity aligns with Bengaluru's irregular fractal spectrum and elevated fractal dimension values, highlighting it as the city with the most pronounced AQI fluctuations.

Mumbai also exhibits a notable disparity (around 82,902), reinforcing its moderate variability and higher AQI volatility compared to Delhi and Kolkata. Kolkata and Delhi show lower disparities of approximately 49,820 and 43,313, respectively, suggesting relatively stable AQI patterns with fewer significant shifts. This consistency between scale disparity results and fractal analyses underscores the utility of examining scale-dependent AQI patterns as an effective approach for assessing risk across urban areas.

## DISCUSSION:

The fractal analysis of AQI values across the four cities allowed for a nuanced understanding of volatility, complexity, and ruggedness in air quality, highlighting patterns of stability and irregularity that traditional statistical methods might not fully capture. However, to draw meaningful conclusions regarding health risks associated with AQI levels, it was necessary to integrate insights from classical statistical metrics alongside fractal analysis results.

**Table 4: Descriptive Statistics of AQI in each city**

| City | N | Mean | Std. Dev. | Minimum | Median | Maximum |
|---|---|---|---|---|---|---|
| Delhi | 2465 | 207.764 | 104.394 | 41 | 193 | 494 |
| Mumbai | 2830 | 106.109 | 54.522 | 25 | 90 | 381 |
| Kolkata | 2375 | 115.764 | 86.038 | 22 | 78 | 419 |
| Bengaluru | 2830 | 73.386 | 25.416 | 24 | 67 | 206 |

### Health Statements for AQI Categories

| AQI | Category | Color Code | Possible Health Impacts |
|---|---|---|---|
| 0-50 | Good | | Minimal Impact |
| 51-100 | Satisfactory | | Minor breathing discomfort to sensitive people |
| 101-200 | Moderate | | Breathing discomfort to the people with lungs, asthma and heart diseases |
| 201-300 | Poor | | Breathing discomfort to most people on prolonged exposure |
| 301-400 | Very Poor | | Respiratory illness on prolonged exposure |
| 401-500 | Severe | | Affects healthy people and seriously impacts those with existing diseases |

**Figure 5: Health Statements for AQI Categories**

From the descriptive statistics (Table 4), we observe notable differences in the mean and standard deviation of AQI values across the cities. Delhi stands out with the highest mean AQI (207.764), indicating a consistently high level of air pollution. According to the AQI health impact categories (Figure 5), this falls into the "Poor" to "Very Poor" range, which means "breathing discomfort to most people on prolonged exposure" or even "respiratory illness on prolonged exposure." The substantial standard deviation (104.394) further points to significant fluctuations, which, combined with fractal analysis showing stable high AQI values, suggests that Delhi poses the highest health risk among the four cities due to both the severity and persistence of its air pollution levels.

Kolkata and Mumbai show moderate mean AQI values (115.764 and 106.109, respectively). This places them in the "Moderate" category, where sensitive individuals might experience breathing discomfort. Both cities exhibit similar volatility patterns, suggesting intermediate health risks. While they experience fluctuations in AQI, the levels do not reach the extremes seen in Delhi, making these cities comparatively less hazardous but still concerning for sensitive populations.

In contrast, Bengaluru has the lowest mean AQI (73.386) and a smaller standard deviation (25.416), placing it in the "Satisfactory" range according to the health impact chart. This suggests generally safer air quality levels, where "minor breathing discomfort" may only affect sensitive individuals. However, fractal analysis reveals high volatility, suggesting that while overall AQI levels are lower, they fluctuate unpredictably. This unpredictability presents challenges for forecasting and policy interventions, as sudden spikes in AQI could lead to temporary health risks despite the overall safer levels. Thus, while Bengaluru is comparatively the safest city for long-term exposure, its volatility requires careful monitoring and adaptive policy measures to mitigate sudden changes in air quality.

Overall, this combined analysis underscores that while Delhi poses the greatest health risks due to consistently high pollution levels, Bengaluru's unpredictable AQI shifts demand agile response strategies from both forecasting agencies and policymakers to ensure effective health risk management.

## CONCLUSION:

This study demonstrates the significant value of fractal analysis in enhancing our understanding of air quality dynamics, providing insights that go beyond traditional statistical approaches. By applying fractal analysis to AQI data from four major Indian cities-Delhi, Mumbai, Kolkata, and Bengaluru-we were able to uncover patterns of volatility and complexity that reveal the underlying behaviours of air pollution, which might otherwise go unnoticed.

While classical statistical metrics helped characterize the general pollution levels across the cities, fractal analysis proved crucial in highlighting scale disparities and the unpredictable fluctuations in AQI. This added depth to our interpretation of health risks, showing that cities like Delhi, with consistently high AQI levels, pose a more persistent health risk, requiring long-term pollution control strategies. In contrast, cities with more volatile AQI, such as Bengaluru, present challenges that require adaptive policies to manage sudden spikes and avoid temporary health risks.

Overall, this combined approach of fractal and statistical analysis provides a more nuanced view of air quality patterns and emphasizes the need for tailored, region-specific interventions. It also highlights the potential of fractal analysis as a valuable tool in environmental research, offering a deeper understanding of complex and variable phenomena like air pollution. This study advocates for the broader application of fractal analysis in air quality management to improve the effectiveness of health risk assessments and policy-making.